\documentclass{gtpart}     % Basic GT/GTM/AGT style
\title{Products of Greek letter elements dug up from the third Morava stabilizer algebra}

\author{Ryo Kato}
\givenname{Ryo}
\surname{Kato}
\address{Graduate school of Mathematics, Nagoya University, Nagoya, 464-8602, Japan}
\email{ryo{\_}kato{\_}1128@yahoo.co.jp}
\author{Katsumi Shimomura}
\givenname{Katsumi}
\surname{Shimomura}
\address{Department of Mathematics, Faculty of Science, Kochi University, 2-5-1, Akebono, Kochi, 780-8520, Japan}
\email{katsumi@kochi-u.ac.jp}

\arxivreference{}
\arxivpassword{}

\volumenumber{}
\issuenumber{}
\publicationyear{}
\papernumber{}
\startpage{}
\endpage{}
\doi{}
\MR{}
\Zbl{}
\received{}
\revised{}
\accepted{}
\published{}
\publishedonline{}
\proposed{}
\seconded{}
\corresponding{}
\editor{}
\version{}

\newtheorem{theorem}{Theorem}[section]
\newtheorem{lemma}[theorem]{Lemma}
\newtheorem{sublemma}[theorem]{Sublemma}
\newtheorem{corollary}[theorem]{Corollary}
\newtheorem{proposition}[theorem]{Proposition}

\theoremstyle{definition}

\theoremstyle{remark}
\newtheorem{remark}[theorem]{Remark}

\newcommand{\lem}[2]{\begin{lemma}\label{#1} #2\end{lemma}}

\newcommand{\thm}[2]{\begin{theorem}\label{#1}#2\end{theorem}}
\newcommand{\cor}[2]{\begin{corollary}\label{#1}#2\end{corollary}}
\newcommand{\rmk}[1]{\begin{remark}#1 \end{remark}}
\newcommand{\pf}[1]{\begin{proof}#1 \end{proof}}
\newcommand{\pfc}[2]{\begin{proof}[Proof of #1]#2\end{proof}}

\renewcommand{\labelenumi}{\arabic{enumi})} %\roman, \Roman, \alph, \Alph

\newcommand{\al}{\alpha}
\newcommand{\be}{\beta}
\newcommand{\ga}{\gamma}

\newcommand{\de}{\delta}
\newcommand{\De}{{\Delta}}

\newcommand{\ze}{\zeta}
\newcommand{\et}{\eta}
\renewcommand{\th}{\theta}

\newcommand{\io}{\iota}
\newcommand{\ka}{\kappa}
\newcommand{\la}{\lambda}

\newcommand{\rh}{\rho}

\newcommand{\Si}{{\Sigma}}

\newcommand{\Om}{{\Omega}}
\newcommand{\barr}{\begin{array}}
\newcommand{\ear}{\end{array}}
\newcommand{\bsk}{\begin{array}{rcl}}
\newcommand{\esk}{\end{array}}
\newcommand{\skh}{\begin{eqnarray*}}
\newcommand{\sko}{\end{eqnarray*}}
\newcommand{\bcss}{\begin{cases}}
\newcommand{\ecss}{\end{cases}}
\newcommand{\sk}[1]{$$\begin{array}{rcl}#1\end{array}$$}
\newcommand{\skr}[2]{$$\begin{array}{rcl}#1\end{array}\lnr{#2}$$}

\newcommand{\AR}[2]{$$\begin{array}{#1}#2\end{array}$$}
\newcommand{\ARr}[3]{$$\begin{array}{#1}#2\end{array}\lnr{#3}$$}

\newcommand{\cass}[1]{\begin{cases}#1\ecss}
\newcommand{\ak}{\quad}

\renewcommand{\Z}{\mathbb Z}
\newcommand{\Zp}{{\mathbb Z}_{(p)}}
\def\o+{\oplus}

\newcommand{\ox}{\otimes}

\newcommand{\e}{{\rm Ext}}
\newcommand{\lnb}{\refstepcounter{theorem}\leqno(\thetheorem)}
\newcommand{\lnr}[1]{\lnb\label{#1}}
\newcommand{\nr}{\refstepcounter{theorem}\thetheorem}
\newcommand{\kko}[1]{(\ref{#1})}
\newcommand{\mbx}[1]{\quad\mbox{#1}\quad}
\newcommand{\cf}{{\it cf.}\ }
\newcommand{\sus}{\subset }
\newcommand{\sm}{\wedge }
\newcommand{\ar}[1]{\xrightarrow{#1} }
\newcommand{\arr}{\to}

\newcommand{\ds}{\displaystyle }
\newcommand{\cg}{\equiv}
\newcommand{\cn}{&\equiv &}

\newcommand{\en}{&= &}
\newcommand{\es}{\ = \ }
\newcommand{\lrk}[1]{\langle #1\rangle}
\newcommand{\qand}{\mbx{and}}
\newcommand{\eR}{\et_R}

\newcommand{\cln}{\colon}
\renewcommand{\O}[1]{\overline{#1}}

\newcommand{\AN}{Adams-Novikov}
\renewcommand{\ss}{spectral sequence}
\newcommand{\ANSS}{Adams-Novikov spectral sequence}

\newcommand{\bet}{\be_{p^2/p^2}}
\newcommand{\bes}{\be_{p/p}}

\newcommand{\xbet}{{^X\!\be}_{p^2/p^2}}

\newcommand{\obes}{\O\be_{p/p}}

\renewcommand{\C}[2]{{#1 \choose #2}}

\newcommand{\p}[1]{^{p^{#1}}
}

\newcommand{\pp}[1]{%^{p^{#1}}
\ifcase #1 1\or p \else {p^{#1}} \fi
}

\newcommand{\q}[1]{
 \count1=1  \advance \count 1 by #1 
\pp{\the\count1}+\pp{#1}}

\begin{document}

\begin{abstract}    % type your abstract below
In \cite{os},
Oka and the second author considered the cohomology of the second Morava stabilizer algebra to study nontriviality of the products of beta elements of the stable homotopy groups of spheres.
In this paper, we use the cohomology of the third Morava stabilizer algebra to find nontrivial products of Greek letters of the stable homotopy groups of spheres:
$\al_1\ga_t$, $\be_2\ga_t$, $\lrk{\al_1,\al_1,\be_{p/p}^p}\ga_t\be_1$  and $\lrk{\be_1,p,\ga_t}$ for $t$ with $p\nmid t(t^2-1)$ for a prime number $p>5$.
\end{abstract}

\maketitle

%%%%%%%%%%%%%%%%%%%%   Start of main body of article

\section{Introduction}

Greek letter elements are well known generators of the stable homotopy groups of spheres localized at a prime $p$. Studying products among these elements is an interesting subject, and studied by several authors.
For example, at an odd prime $p$,  all products of alpha elements are trivial.
In \cite{os},
we used $H^*S(2)$ to study nontriviality of the product of beta elements.
In this paper, we use $H^*S(3)$ to find relations of Greek letters.
The multiplicative structure of $H^*S(3)$ is given by Yamaguchi \cite{y}, but unfortunately, it has some typos. %Indeed,
%his table of relations implies a contradiction $-3h_1a_2b_0=g_0'a_0a_2=-2h_1a_2b_0$.
So here, our computation is based on Ravenel's.

Let $\bes$ be the generator of the $E_2$-term $E_2^{2,p^2q}(S)$ of the \ANSS\ converging to the homotopy groups $\pi_*(S)$ of the sphere spectrum $S$.
Hereafter, $q=2p-2$ as usual.
A relation given by Toda implies that $\bes$ dies in the \ANSS\ at a prime $p>2$.
At the prime two, $\be_{2/2}^2=0$ by \cite[Prop. 8.22]{mrw}, while
at the prime numbers three and five, Ravenel showed that $\be_{p/p}^p$ survives to a homotopy element of $\pi_*(S)$ and $\al_1\bes^p=0$ for the generator $\al_1$ of $\pi_{q-1}(S)$.
Here, we show the following

\thm{bes}{At a prime $p> 3$, $\bes^p$ survives to $\pi_{(p^3-1)q-2}(S)$ and $\al_1\bes^p=0$. }

\cor{bes:1}{At a prime $p> 3$, the Toda bracket $\lrk{\al_1,\al_1,\bes^p}\left(=\al_1\beta_{p^2/p^2}\right)$ is defined.}

\rmk{
It is already known that $\al_1\beta_{p^2/p^2}$ survives in the Adams-Novikov spectral sequence by the work of R. Cohen \cite{c}.
Corollary \ref{bes:1} states that the Cohen's element is a Toda bracket $\lrk{\al_1,\al_1,\bes^p}$.
}

We notice that at the prime $3$, Ravenel showed these in \cite{r:book}.

Let  $\be_1$, $\be_2$ and $\ga_t$ $(t>0)$ be the generators of {Coker~\!$J$} of dimensions $pq-2$, $(2p+1)q-2$ and $(tp^2+(t-1)p+t-2)q-3$, respectively. 

\thm{main}
{Let $p>5$, and $t$ be a positive integer with $p\nmid t(t^2-1)$.
Then, the elements $\al_1\ga_t$, $\be_2\ga_t$, $\lrk{\al_1,\al_1,\be_{p/p}^p}\be_1\ga_t$  and $\lrk{\be_1,p,\ga_t}$ generate  subgroups of the stable homotopy groups of spheres isomorphic to $\Z/p$.
Besides, even in the case $p|(t+1)$, $\be_1\ga_t$ and $\lrk{\be_1,p,\ga_t}$ are generators of order $p$.}

Note that $\lrk{\be_1,p,\ga_t}=\lrk{\ga_t,p,\be_1}$.
%This result contains the first $\de$-element $\de_1=\al_1\ga_{p-1}$.
We also notice that if $t=1$, then $\lrk{\ga_1,p,\be_1}=0$, while $\be_2\ga_1$ is non-trivial (see section five).

From here on, we assume that the prime number $p$ is greater than three.

\section{$H^*S(3)$ revisited}

We begin with recalling some notation from Ravenel's green book \cite{r:book}.
Let $BP$ denote the Brown-Peterson spectrum. Then, the pair
$$
	(BP_*,BP_*(BP))=(\Zp[v_1,v_2,\dots], BP_*[t_1,t_2,\dots])
$$
is a Hopf algebroid. Here, the degrees of $v_i$ and $t_i$ are $2p^i-2$.
The structure maps act as follows:
\skr{
	\eR(v_1)\en v_1+pt_1\\
	\eR(v_2)\cn v_2+v_1t_1^p+pt_2\mod (p^2,v_1^p)\\
	\eR(v_3)\cn v_3+v_2t_1\p 2+v_1t_2^p+pt_3\\
&&\hspace{.3in}-pv_1v_2^{p-1}(t_2+t_1^{p+1})\mod (p^2,v_1^2,v_2^p)\\
	\De(t_1)\en t_1\ox 1+1\ox t_1\\
	\De(t_2)\en t_2\ox 1+t_1\ox t_1^p+1\ox t_2-v_1b_{10}\\
	\De(t_3)\cn t_3\ox 1+t_2\ox t_1\p2+t_1\ox t_2^p+1\ox t_3\mod (v_1,v_2)\\
	\De(t_4)\cn t_4\ox 1+t_3\ox t_1\p3+t_2\ox t_2\p2+t_1\ox t_3^p+1\ox t_4\\
&&\hspace{.3in}-v_3b_{12}\mod (v_1,v_2)\\
	\Delta(t_5)&=&t_5\ox 1+t_4\ox t_1^{p^4}+t_3\ox t_2^{p^3}+t_2\ox t_3^{p^2}+t_1\ox t_4^p+1\ox t_5\\
&&\hspace{.3in}-v_3b_{22}-v_4b_{13} \mod (p,v_1,v_2)
}{str}
for 
\ARr{rcl}{
	b_{1k}\en \dfrac 1p\left(\Delta(t_1)^{p^{k+1}}-t_1^{p^{k+1}}\otimes 1-1\ox t_1^{p^{k+1}}\right)\es \ds \frac1p\sum_{i=1}^{p^{k+1}-1}\C {p^{k+1}} it_1^i\ox t_1^{p^{k+1}-i}
\qand\\
b_{2k}\en \dfrac 1p\left(\Delta(t_2)^{p^{k+1}}-t_2^{p^{k+1}}\otimes 1-t_1^{p^{k+1}}\ox t_1^{p^{k+2}}-1\ox t_2^{p^{k+1}}-v_1^{p^{k+1}}b_{1k+1}\right).
}{bi}

Let $K(3)_*=F_p[v_3,v_3^{-1}]$ have the $BP_*$-module structure given by $v_iv_3^s=v_3^sv_i=v_3^{s+1}$ if $i= 3$, and $=0$ otherwise, and 
\sk{
	\Si(3)\en K(3)_*\ox_{BP_*}BP_*(BP)\ox_{BP_*}K(3)_*\\
	\en K(3)_*[t_1,t_2,\dots]/(v_3t_i\p3-v_3\p it_i:i>0)\mbx{(by \cite[6.1.16]{r:book})}
}
is the Hopf algebra with structure inherited from $BP_*(BP)$.
Define the Hopf algebra $S(3)$ by $S(3)=\Si (3)\ox _{K(3)_*}\!\! F_p$, where $K(3)_*$ acts on $F_p$ by $v_3\cdot 1=1$. Then, 
$$
	S(3)=F_p[t_1,t_2,\dots]/(t_i\p3-t_i:i>0).
$$
Now we abbreviate $\e_{S(3)}(F_p,F_p)$ to $H^*S(3)$. 

Consider integers $d_i$ ($=d_{3,i}$ in \cite[6.3.1]{r:book})
$$
	d_{i}=\cass{0& i\le 0,\\ \max(i,pd_{i-3})&i>0.}
$$
Then, there is a unique increasing filtration of the Hopf algebroid $S(3)$ with deg $t_i^{p^j}=d_{i}$ for $0\le j<3$.

\thm{r:632}{{\rm (Ravenel\cite[6.3.2]{r:book})} The associated Hopf algebra $E^0S(3)$ is isomorphic to the truncated polynomial algebra of height $p$ on the elements $t_i^{p^j}$ for $i>0$ and $j\in \Z/3$, with coproduct defined by
$$
	\De (t_i^{p^j})=\cass{\sum_{k=0}^i t_k^{p^j}\ox t_{i-k}^{p^{k+j}}&i\le 3,\\
		t_i^{p^j}\ox 1+1\ox t_i^{p^j}+b_{i-3,j+2}&i>3.}
$$}

Let $L(3)$ be the Lie algebra without restriction with basis $x_{i,j}$ for $i>0$ and $j\in \Z/3$ and bracket given by 
$$
	[x_{i,j},x_{k,l}]=\cass{\de_{i+j}^lx_{i+k,j}-\de_{k+l}^jx_{i+k,l}&\mbox{for $i+k\le 3$,}\\0&\mbox{otherwise,}} %\ak \LR{\de_j^i=\cass{1&i=j\\0&i\ne j}}
$$
where $\de_j^i=1$ if $i\cg j$ mod $3$ and $0$ otherwise,
%\footnote{In the green book, the second $\de$ is $\de_{k+1}^j$, which is a typo%. Indeed, $\lrk{x_{ij}x_{kl}, t_{mn}}=\lrk{x_{ij}\ox x_{kl}, \De (t_{mn})}=\lrk{x_{ij}\ox x_{kl}, \sum_{a=0}^m t_{an}\ox t_{m-a,a+n}}$ and so $i=a$, $j=n$, $k=m-a$ and $l=a+n$, which shows that $m=i+k$, $n=j$ and $l=i+j$. This gives the first term. By switching $(ij)$ and $(kl)$, the second follows from $m=i+k$, $n=l$ and $j=k+l$.}
and $L(3,k)$ the quotient of $L(3)$ obtained by setting $x_{i,j}=0$ for $i>k$.
Then, Ravenel noticed in \cite[6.3.8]{r:book}:

\thm{638}{$H^*(L(3,k))$ for $k\le 3$ is the cohomology of the exterior complex $E(h_{i,j})$on one-dimensional generators $h_{i,j}$ with $i\le k$ and $j\in \Z/3$, with coboundary
$$
	d(h_{i,j})=\sum_{s=1}^{i-1}h_{s,j}h_{i-s,s+j}.
$$
}
From now on, we abbreviate $h_{i,j}$ to $h_{ij}$, and $h_{1j}$ to $h_j$.

Under the above filtration, Ravenel constructed the May spectral sequences

\thm{Ravenel-May}{{\rm (Ravenel \cite[6.3.4, 6.3.5]{r:book})} There are spectral sequences
\renewcommand{\labelenumi}{(\alph{enumi})}
\begin{enumerate}
\item $E_2=H^*(L(3,3))\Longrightarrow H^*(E_0S(3))$ and
\item $E_2=H^*(E_0S(3))\Longrightarrow H^*(S(3))$.
\end{enumerate}
}
Since these spectral sequences collapse, $H^*S(3)$ is additively isomorphic to $H^*L(3,3)$.
Therefore, we have a projection 
$$
	\pi\cln H^*S(3)\to E^0H^*S(3)=H^*(E_0S(3))=H^*L(3,3).
\lnr{proj}$$
Note that the Massey product $\lrk{h_i, h_{i+1}, h_{i+2},h_i}$ is homologous to $v_3^{(2-p)p^i}b_{i+2}$ represented by $v_3^{(2-p)p^i}b_{1,i+2}$ of \kko{bi}, and $\pi$ assigns the Massey product to $b_{i+2}\in H^*L(3,3)$.
Ravenel determined in \cite[6.3.34]{r:book} the additive structure of $H^*L(3,3)$.
In particular, we have the following:

\thm{r:6334}{ $H^*L(3,3)$ contains submodules generated by the linear independent elements:
$$
 h_1k_1\ze_3,\ak b_0k_1\ze_3,\ak h_0l, \ak k_0l, \ak h_0b_0b_2l \qand h_1l.
$$
%Moreover $h_1l\ne h_1k_1\ze_3$.
Here, $l=h_2h_{21}h_{30}$, 
%$k_i=h_{i+1}h_{2i}$ $(i=0,1)$,
$k_i=h_{2i}h_{i+1}$ $(i=0,1)$,
 $b_0=h_1h_{32}+h_{21}h_{20}+h_{31}h_1$, $b_2=h_0h_{31}+h_{20}h_{22}+h_{30}h_0$ and $\ze_3=h_{30}+h_{31}+h_{32}$.
}

\pf{
In the table of the proof of \cite[6.3.34]{r:book}, we find the elements 
$$
	h_0,\ak h_1,\ak k_0,\ak  b_0,\ak  b_2,\ak  l, \ak l'=h_0h_{22}h_{31},%\ak h_1k_1h_{31}
\qand \ze_3,
$$
as well as the first element $ h_1k_1\ze_3$ of the theorem.
We also have the element $-h_1k_1h_{30}=h_1h_2h_{21}h_{30}$ in the table, which is the last element $h_1l$ of the theorem.
%These also imply $h_1l\ne h_1k_1\ze_3$.
Besides $h_1k_1h_{31}$ and $h_1k_1h_{32}$ are in the table too. 
We see that %$b_0k_1=h_1k_1\ze_3-h_1k_1h_{30}$, 
$b_0k_1=-h_1k_1h_{31}+h_1k_1h_{32}$
and so the second
 element is given by  $b_0k_1\zeta_3=-h_1k_1h_{31}\ze_3+h_1k_1h_{32}\ze_3$.% given by
%\[
%\begin{array}{rcl}
%b_0k_1\zeta_3 &=& %(h_1h_{32}+h_{31}h_1)h_{21}h_2\zeta_3 \\
%h_1k_1(h_{32}-h_{31})\zeta_3 
%h_1l\zeta_3 -2h_1k_1h_{31}\zeta_3 .
%\end{array}
%\]

The element $h_0b_0b_2l\ze_3$ is computed as
\sk{
	&& h_0h_2h_{21}h_{30}(h_1h_{32}+h_{21}h_{20}+h_{31}h_1)(h_0h_{31}+h_{20}h_{22}+h_{30}h_0)(h_{30}+h_{31}+h_{32})\\
	%\en h_0h_1h_2h_{21}h_{20}h_{22}h_{30}(h_{32}-h_{31})(h_{31}+h_{32})\\
	%\en h_0h_1h_2h_{21}h_{20}h_{22}h_{30}(h_{32}h_{31}-h_{31}h_{32})\\
	\en -2h_0h_1h_2h_{20}h_{21}h_{22}h_{30}h_{31}h_{32}.
}
Therefore, $h_0b_0b_2l$ is the dual of the generator $-\frac12\ze_3$, and the elements $h_0b_0b_2l$ and $h_0l$ are generators.
%Put $l'=h_0h_{22}h_{31}$, which is the generator in the table.
Similarly, a computation
\sk{
	k_0ll'\ze_3\en h_{20}h_{1}h_2h_{21}h_{30}h_0h_{22}h_{31}(h_{30}+h_{31}+h_{32})\\
	\en h_0h_1h_2h_{20}h_{21}h_{22}h_{30}h_{31}h_{32}\\
}
shows that $k_0l$ is the dual of the generator $l'\ze_3$.
}

\lem{trivial}{In $H^*L(3,3)$, $h_0k_1=0$ and $k_0k_1=0$. }

\pf{From the proof of \cite[6.3.34]{r:book}, we read off the relations
$h_0k_1=e_{30}h_2$ and $k_0k_1=e_{30}g_1$ in $H^*L(3,2)$.
Since $e_{30}$ cobounds $h_{30}$ in $H^*L(3,3)$, the lemma follows.
}

\section{Greek letter elements}

Let $E_r^{s,t}(X)$ denote the $E_r$-term of the \ANSS\ converging to the homotopy group $\pi_{t-s}(X)$ of a spectrum $X$.
Then the $E_2$-term is $\e_{BP_*(BP)}(BP_*,BP_*(X))$.
We here consider the Ext-group $\e_{BP_*(BP)}(BP_*,M)$ for a $BP_*(BP)$-comodule $M$ as the cohomology of the cobar complex $\Om_{BP_*(BP)}^*M$ (\cf \cite{mrw}).
Consider a sequence $A=(a_0, a_1, \dots, a_n)$ of non-negative integers so that  the sequence $p^{a_0},v_1^{a_1},\dots, v_n^{a_n}$ is invariant and regular.
For such a sequence $A$, 
 Miller, Ravenel and Wilson introduced in \cite{mrw} $n$-th Greek letter elements $\et^{(n)}_{s(A)}$ in the \AN\ $E_2$-term $E_2^{n,t(A)}(S)$ by
\skr{
	\et^{(n)}_{s(A)}\en \de_{A,1}\cdots\de_{A,n}(v_n^{a_n})\in E_2^{n,t(A)}(S)
}{def:G}
for $v_n^{a_n}\in \e_{BP_*(BP)}^{0,2a_n(p^n-1)}(BP_*,BP_*/I(A,n))$.
Here, $s(A)=a_n/a_{n-1},a_{n-2},\cdots, a_0$ and $t(A)=2a_n(p^n-1)-2\sum_{i=0}^{n-1}a_i(p^i-1)$, $I(A,k)$ denotes the ideal of $BP_*$ generated by $p^{a_0},v_1^{a_1},\dots, v_{k-1}^{a_{k-1}}$, and $\de_{A,k+1}$ is
 the connecting homomorphism  associated to the short exact sequence
$$
	0\arr BP_*/I(A,k)\ar{v_k^{a_k}}BP_*/I(A,k)\arr BP_*/I(A,k+1)\arr 0.
$$
In particular, we write $\al=\et^{(1)}$, $\be=\et^{(2)}$ and $\ga=\et^{(3)}$.
So far, only when $n\le 3$, 
%we know a condition whether or not Greek letter elements survive to homotopy elements.
many conditions for that Greek letter elements survives to homotopy elements are known.
We abbreviate $\et^{(n)}_{s(A)}$ to $\et^{(n)}_{a_n}$ if $A=(1,\dots,1,a_n)$ as usual.
For example, we consider $\be$-elements defined by 
$$\begin{array}{l}
	\be_s=\de_{(1,1),1}(\be_s')\in E_2^{2,t(1,1,s)}(S)\\
	\qquad\mbx{for $\be_s'=\de_{(1,1),2}(v_2^s) \in E_2^{1,t(1,1,s)}(V(0))$, and}\\
	\be_{p^i/p^i}=\be_{p^i/p^i,1}=\de_{(1,p^i),1}\de_{(1,p^i),2}(v_2^{p^i})
	\in E_2^{2,t(1,p^i,p^i)}(S).
\end{array}\lnr{beta}$$

%At the prime $p$ greater than three, 
Hereafter we assume that the prime $p$ is greater than three.
We have the Smith-Toda spectrum $V(k)$ for $k=0,1,2$ defined by the cofiber sequences
\ARr{c}{
	S\ar{p} S\ar{i} V(0)\ar{j} \Si S,\\[1mm]
	\Si ^{q}V(0)\ar{\al} V(0)\ar{i_1} V(1)\ar{j_1} \Si ^{q+1}V(0)\qand \\[1mm]
	\Si ^{(p+1)q}V(1)\ar{\be} V(1)\ar{i_2} V(2)\ar{j_2} \Si ^{(p+1)q+1}V(1).
}{cof}
Here, $\al\in [V(0),V(0)]_q$ is the Adams map and $\be\in [V(1), V(1)]_{(p+1)q}$ is the $v_2$-periodic element due to L. Smith.
Note that the $BP_*$-homology of these spectra are $BP_*(V(k))=BP_*/I_{k+1}$ for the ideal $I_k$ of $BP_*$ generated by $v_i$ for $0\le i<k$ with $v_0=p$.
We consider the Bousfield-Ravenel localization functor $L_3$ with respect to $v_3^{-1}BP$.
%Let $S(3)$ denote the dual of the Morava stabilizer algebra.
The $E_2$-term  $E_2^*(L_3V(2))$ of $L_3V(2)$ is isomorphic to $K(3)_*\ox H^*S(3)$, whose structure is given in \cite{r:book} (see also \cite{y}),
%Here, $K(3)_*=F_p[v_3^{\pm 1}]$.
and we consider the composite
$$
	r\cln E_2^*(S)\ar{\io_*} E_2^*(V(2))\ar{\et} E_2^*(L_3V(2))\ar{\rh} H^*(S(3))\ar{\pi} H^*L(3,3).
$$
Here the first map is induced from the inclusion $\io\cln S\to V(2)$ to the bottom cell, the second is from the localization map, the third is obtained by setting $v_3=1$ and the last is the projection \kko{proj}.

\lem{rG}{The map $r$ assigns the Greek letter elements as follows:
\sk{
	r(\al_1)\en h_0,
%	r(\be_1')\en h_1, 
	\ak r(\be_1)\es -b_0,\ak
%	r(\be_2)\es -2k_0,\\
	r(\be_2)\es 2k_0,\\
%	r(\ga_t')\en t(t-1)k_1, 
%	r(\ga_t)\en t(t^2-1)l-t(t-1)k_1\ze_3\qand
	r(\ga_t)\en -t(t^2-1)l-t(t-1)k_1\ze_3\qand
	r(\bes)\es -b_1.
}
%The other elements are assigned to zero.
We also have $\be'_1=h_1-v_1^{p-1}h_0\in E_2^{1,pq}(V(0))$ for the generators $h_i$ of $E_2^{1,p^iq}(V(0))$ represented by $t_1^{p^i}$.
}
%Here,
%\sk{
%	h_i\en [t_1\p i], \ak b_0\es [t_1^p\ox t_3\p2+t_2^p\ox t_2+t_3^p\ox t_1^p]\\
%	k_i\en [t_2\p i\ox t_1\p{i+1}+\frac12 t_1\p i\ox t_1^{2p^{i+1}}]\\
%	l\en [t_1\p2\ox t_2^p\ox t_3+\dots]\in\lrk{k_1,h_2,h_1,h_0}.
%}

\pf{
	First we consider the images of the Greek letter elements under the map $\io_*\cln E_2^*(S)\to E_2^*(V(2))$.
In the cobar complex $\Om_{BP_*(BP)}^*BP_*$, by \kko{str},
$d(v_1)= pt_1$, 
%$d(v_2\p i)\cg v_1\p it_1\p{i+1}$ mod $(p,v_1\p{i+1})$
$d(v_2\p i)\cg v_1\p it_1\p{i+1}-v_1^{p^{i+1}}t_1^{p^i}$ mod $(p)$
 for $i\ge 0$,  $d(v_2^2)\cg 2v_1v_2t_1^p+v_1^2t_1^{2p}$ mod $(p,v_1^p)$,  and $d(v_3^t)\cg tv_2v_3^{t-1}t_1\p2+\C t2 v_2^2v_3^{t-2}t_1^{2p^2}+{t\choose 3}v_2^3v_3^{t-3}t_1^{3p^2}$ mod $(p,v_1, v_2^4)$, which imply
\sk{
	\de_{(1),1}(v_1)\en [t_1], \ak \de_{(1,1),2}(v_2)\es [t_1^p-v_1^{p-1}t_1], \\[1mm]
	\de_{(1,1),2}(v_2^2)\en [2v_2t_1^p+v_1t_1^{2p}+v_1^{p-1}y],%\es 2v_2h_1\\
	\ak \de_{(1,p),2}(v_2^p)\es [t_1\p2-v_1^{p^2-p}t_1^p]\qand\\ \de_{(1,1,1),3}(v_3^t)\en [tv_3^{t-1}t_1\p2+\C t2 v_2v_3^{t-2}t_1^{2p^2}+\C t3v_2^2v_3^{t-3}t_1^{3p^2}+v_2^3z]\es \O\ga_{t},
}
for cochains $y\in \Om_{BP_*(BP)}^1BP_*/(p)$ and $z\in\Om_{BP_*(BP)}^1BP_*/(p,v_1)$.
Here, $[x]$ denotes a cohomology class represented by a cocycle $x$.
The first one shows $\al_1=h_0$, and the second gives the last statement of the lemma.
We further see that $d(t_1\p k)=-pb_{1k-1}$ for $k\ge 1$ and $d(v_k)\cg pt_k$ mod $I((2,1,1),k)$ for $k=2,3$ by \kko{str} in $\in \Om_{BP_*(BP)}^1BP_*$.
Moreover, $[b_{1k}]$'s are assigned to $b_k$ in $H^*L(3,3)$ under the projection $\pi$,  and we obtain 
\AR{l}{
	r\de_{(1,p^{k-1}),1}(h_k-v_1^{p^k-p^{k-1}}h_{k-1})\es -b_{k-1} \mbx{for $k=1,2$,}\\ r\de_{(1,1),1}([2v_2t_1^p+v_1t_1^{2p}])% \es r[2t_2\ox t_1^p+2v_2b_0+t_1\ox t_1^{2p}+\dots]
	\es 2k_0,\\ %\ak r\de_{(1,p),1}(h_2)\es -b_1, \\
  \de_{(1,1,1),2}(\O\ga_{t})\es [t(t-1)v_3^{t-2}t_2^p\ox t_1\p2+{t\choose 2}v_3^{t-2}t_1^p\otimes t_1^{2p^2}+w]\es \ga'_t \qand \\[1mm]
%	 \es t(t-1)v_3^{t-3}[t_2^p\ox v_3t_1\p2+\frac12 v_3t_1^p\ox t_1^{2p^2}+v_1v_2z'+v_2^2z'']\es \ga'_t \qand \\[1mm]
	r\de_{(1,1,1),1}(\ga'_t)\es t(t-1)(t-2)h_{30}k_1+t(t-1)r\de_{(1,1,1),1}([t_2^p\otimes t_1^{p^2}+\frac 12t_1^p\otimes t_1^{2p^2}]).
}
Here, $w$ is a linear combination of terms in the ideal $(v_1,v_2)^2$.
% and of the form $v_ex\ox y$ for $e\in\{1,2\}$ and $x,y\in K(3)_*\{t_i\p kt_j\p l, t_1^{3p^2}: i,j,k,l\in\{1,2\}\}$.
Thus the relations other than $r(\ga_t)$ follows.
%Note that $b_1=h_2h_{30}+h_{22}h_{21}+h_{32}h_2$. %'±'±'Í•K—v'Å'µ'傤'©?

We note that $b_{20}$ in \kko{bi} corresponds to $h_{21}h_{30}+h_{31}h_{21}$ by $\Delta(t_5)^p$ in \kko{str}.
Since $d(t_2^p)=-t_1^p\ox t_1^{p^2}+v_1^pb_{11}-pb_{20}$ by (\ref{str}),
we obtain $r\de_{(1,1,1),1}([t_2^p\otimes t_1^{p^2}+\frac 12t_1^p\otimes t_1^{2p^2}])=-(h_{21}h_{30}+h_{31}h_{21})h_2+h_{21}b_1=-3l-k_1\ze_3$, which shows
%by $\Delta(t_5)^p$ and $d(t_2^p)=-t_1^p\ox t_1^{p^2}+v_1^pb_{11}-pb_{20}$ (see (\ref{str})), we obtain 
the relation on $r(\ga_t)$.
}

Recall the cofiber sequences \kko{cof} and the $v_3$-periodic element $\ga\in [V(2),V(2)]_{q_3}$ ($q_3=(p^2+p+1)q$) due to H. Toda.
Then, the Greek letter elements in homotopy are defined by
\skr{
	\al_t\en j\al^ti,\ak 
	\be_t\es j\be_t' \mbx{for $\be'_t\es j_1\be^ti_1i$\ak and }
	\ga_t\es jj_1j_2\ga^ti_2i_1i}{Greek}
for $t>0$, and the Greek elements in the $E_2$-term survives to the same named one in homotopy by the Geometric Boundary Theorem (\cf \cite{r:book}).

\pfc{Theorem \ref{main}}{
We begin with noticing that the element $b_i$ in $H^*L(3,3)$ is the image of the Massey product $\lrk{h_i,h_{i+1}, h_{i+2}, h_i}$ under $\pi$, which is homologous to $b_i$ represented by $b_{1i}$ in \kko{bi}.
We further note that the Toda brackets $\lrk{\al_1,\al_1,\bes^p}$ and $\lrk{\be_1,p,\ga_t}$ are detected by $\al_1b_{2}$ and $h_1\ga_t$ of $E_2^*(S)$, respectively.
Indeed, in the first bracket, $d_{2p-1}(b_2)=\al_1\bes^p$ by Corollary \ref{d2p-1} below, and in the second bracket, $\lrk{\be_1,p,\ga_t}=j\lrk{\be_1',p,\ga_t}$.
Under the condition on $t$,
Lemmas \ref{rG}, \ref{r:6334} and \ref{trivial} imply that each element of $\al_1\ga_t$, $\be_2\ga_t$, $\al_1b_2\ga_t\be_1$ and $h_1\ga_t$, as well as $\be_1\ga_t$, generates a submodule isomorphic to $\Z/p$  of the $E_2$-term $E_2^*(S)$.
These are, of course, permanent cycles, and nothing kills them in the \ANSS\
since each element has a filtration degree less than $2p-1$.
}

\section{$\bes^p$ in the homotopy of spheres}

Let $X$ and $\O X$ be the $(p-1)q$- and $(p-2)q$-skeletons of the Brown-Peterson spectrum $BP$. Then, we have  the cofiber sequences
$$
	S \ar{\io} X\ar{\ka} \Si^q \O X\ar{\la} S^1\qand \O X\ar{\io'}X\ar{\ka'} S^{(p-1)q}\ar{\la'} \Si \O X.
\lnr{cof12}$$

Then, 
$$
	BP_*(X)=BP_*[x]/(x^p) \qand BP_*(\O X)=BP_*[x]/(x^{p-1}) 
$$
as subcomodules of $BP_*(BP)$, where $x$ corresponds to $t_1$.
From \cite[Chap.7]{r:book}, we read off the following:\\

\noindent
(\nr\label{vanish:X})\quad {\it %\begin{minipage}{5in}
%\begin{enumerate}
%\item the vanishing line 
%$$
%	E_2^{2s+e,tq}(X)=0 \mbx{if $t<sp^2+ep$}
%$$
%for $s\ge 0$ and $e=0,1$.
%\item 
		$b_1^p=0\in E_2^{2p, p^3q}(X)$, which implies
$$
	E_2^{2s+e,tq}(X)=0 \mbx{if $s\ge p$ and $t<(s-1)p^2+(s+1+e)p$.}
$$
%\end{enumerate}
%\end{minipage}\\
}
\lem{b0}{$b_0\cln E_2^{2s+e,tq}(S)\to E_2^{2s+2+e,(t+p)q}(S)$ is monomorphic if  $s\ge p$ and $t\le (s-1)p^2+(s+e)p$.}

\pf{
	Note that $b_0=\la\la'$, and the lemma follows from \kko{vanish:X} and the exact sequences
\AR{c}{
	E_2^{2s+e,(t+p-1)q}(X)\ar{\ka'} E_2^{2s+e,tq}(S)\ar{\la'}E_2^{2s+1+e,(t+p-1)q}(\O X)\\%\ar{\io'} E_2^{s+1,t+(p-1)q}(X)\\
	E_2^{2s+e+1,(t+p)q}(X)\arr E_2^{2s+e+1,(t+p-1)q}(\O X)\ar{\la}E_2^{2s+2+e,(t+p)q}(S)%\arr E_2^{s+2,t+pq}(X)
}
induced from the cofiber sequences in \kko{cof12}.
}

%\pf{
%	 $E_2^{2p,p^3q}(X)$ is a subquotient of $\Z/p\{b_1^p\}$, and we compute
%\sk{
%	b_1^p\en b_1^{p-1}\lrk{h_1,\dots,h_1} \es \lrk{b_1^{p-1}, h_1,\dots,h_1}h_1 \mbx{($p$-fold Massey product)}
%}
%in which $\lrk{b_1^{p-1}, h_1,\dots,h_1}$ is killed by $b_{20}^{p-1}$.
%}

Ravenel showed that $d_{2p-1}(\bet)\cg \al_1\bes^p$ mod Ker $\be_1^p$ in the \ANSS\ for $\pi_*(S)$ (\cf \cite[6.4.1]{r:book}).
Here, the mapping $\be_1^p$ on $E_2^{2p+1,(p^3+1)q}(S)$ %$\be_1^p\cln E_2^{2p+1,(p^3+1)q}(S)\to E_2^{4p+1,(p^3+p^2+1)q}(S)$
 is a monomorphism by Lemma \ref{b0}:

\cor{d2p-1}{In the \ANSS\ for $\pi_*(S)$, $d_{2p-1}(\bet)= \al_1\bes^p\in E_{2p-1}^{2p+1,(p^3+1)q}(S)=E_2^{2p+1,(p^3+1)q}(S)$.}

\pfc{Theorem \ref{bes}}{
Consider the first cofiber sequence in \kko{cof12}.
Since the \AN\ $E_2$-term $E_2^{sq+3 ,(p^3+s)q}(X)$ vanishes for $s>0$ by \kko{vanish:X},
 the element $\io_*(\bet)$ $\in E_2^{2,p^3q}(X)$ survives to a homotopy element $\xbet\in \pi_*(X)$.
In general, we see that\\

\noindent
(\nr\label{al})\quad {\it
Let $\O\io\cln S\to \O X$ denote the inclusion to the bottom cell.
Then, $\la_*\O \io (x)=\al_1x$ for $x\in E_2^*(S)$.}
\\

\noindent
Put $\obes=\O\io_*(\bes)\in E_2^{2,p^2q}(\O X)$, and we see that $\la_*(\obes^p)=\al_1\bes^p$, and so
 we see that $\obes^p$ detects an essential homotopy element $\ka_*(\xbet)\in \pi_*(\O X)$ by Corollary \ref{d2p-1} and \cite{s}, which
 we also denote by $\obes^p$.

Now turn to the second cofiber sequence in \kko{cof12}.
%We also see that $\bes^p=0\in E_2^{2p, p^3q}(X)$ by the computation
%\sk{
%	b_1^p\en b_1^{p-1}\lrk{h_1,\dots,h_1} \es \lrk{b_1^{p-1}, h_1,\dots,h_1}h_1 \mbx{($p$-fold Massey product)}
%}
%in which $\lrk{b_1^{p-1}, h_1,\dots,h_1}$ is killed by $b_{20}^{p-1}$.
%By \kko{vanish:X}, we see that $\io'_*(\obes^p)=0\in \pi_*(X)$, and so there is an element $\xi\in \pi_*(S)$ such that $\la'(\xi)=\obes^p$.
%Then, $\be_1\xi=\al_1\bes^p$.
The relation $b_1^p=0$ of \kko{vanish:X} yields a cochain
$y=\sum_{i=0}^{p-1} x^iy_i\in \Om^{2p-1}BP_*(X)$ such that $d(y)=b_1^p$, where $y_i\in \Om^{2p-1} BP_*$.
It follows that $d(\O y)=b_1^p-d(x^{p-1})y_{p-1}\in \Om^{2p}BP_*(\O X)$ for $\O y=\sum_{i=0}^{p-2} x^iy_i\in \Om^{2p-1}BP_*(\O X)$.
In particular $d(y_{p-1})=0\in \Om^{2p-1}BP_*$ and $d(y_{p-2})=(1-p)t_1\ox y_{p-1}$. 
By definition, these imply $\la'_*(y_{p-1})=b_1^p$.
Consider the exact sequence obtained by applying the homotopy groups to the second cofiber sequence.
Then, $\io'_*(\obes^p)=0$ by \kko{vanish:X}, and so $\obes^p$ must be pulled back to an element $\xi\in \pi_*(S)$ detected by $y_{p-1}$.
Since $b_0=\la\la'$, $b_0y_{p-1}=h_0b_1^p$, and 
$
	\lrk{h_0,\dots,h_0}y_{p-1}=h_0\lrk{h_0,\dots,h_0,y_{p-1}}
$,
we see that
$$
	b_1^p\cg \lrk{h_0,\dots,h_0,y_{p-1}}\not\cg 0\in E_2^{2p,p^3q}(S) \mod \ker h_0.
$$
Put $b_1^p= \lrk{h_0,\dots,h_0,y_{p-1}}+c$ for $c\in \ker h_0\sus E_2^{2p, p^3q}(S)$.
Then, $b_1^p-c$ survives to $\bes^p\in \pi_*(S)$.
%Then, $\la_*\O \io (c)=h_0c=0$ for $\la_*\cln E_2^{2p, p^3q}(\O X)\to E_2^{2p+1, (p^3+1)q}(S)$, which is a monomorphism by \kko{vanish:X}.
%It follows that $\O\io(c)=0$.
%So we have an element $w=\sum_{i=0}^{p-3}x^iw_i$ such that $\de(w)=c$.
%$d(\wt w)=\sum_{i=0}^{p-3}\frac1{i+1}d(x^{i+1}w_i)=c$.
%\sk{
%	d(w_{p-3})\en 0\\
%	d(w_{p-4})\en -t_1\ox w_{p-3}\\
%	d(w_{p-5})\en \frac1{p-2}\C {p-2}2t_1^2\ox w_{p-3}+t_1\ox w_{p-4}\\
%	c\en \sum_{i=0}^{p-3}\frac1{i+1}t_1^{i+1}\ox w_i
%}
%Therefore, $c=\langle \overbrace{h_0,\dots,h_0}^{p-2},w_{p-3}\rangle$ for $w_{p-3}\in E_2^{2p-1, (p^3-p+2)q}(S)$.
%
% we have $\wt c\in E_2^{2p,p^3q}(X)$ such that $\ka_*(\wt c)=c$.
%It follows that $\bes^p=\lrk{\al_1,\dots,\al_1,\xi}$ is a homotopy element.
%}

%\lem{al1}{$\al_1\bes^p=0$ in $\pi_*(S)$.}
%
%\pf{
The element
$\al_1\bes^p$ is detected by $h_0(b_1^p-c)=h_0b_1^p$ in the \AN\ $E_2$-term, which is killed by $b_2$ by Corollary \ref{d2p-1}.
}
\section{Remarks}
\subsection{A relation on Toda bracket}
The relation $\lrk{\be_s,p,\ga_t}=\lrk{\ga_t,p,\be_s}$ follows immediately from  results of Toda:
By definition, $\lrk{\be_s,p,\ga_t}=j\be_{(s)}\ga_{(t)}i$ and $\lrk{\ga_t,p,\be_s}=j\ga_{(t)}\be_{(s)}i$
for $\be_{(s)}=j_1\be^si_1$ and $\ga_{(t)}=j_1j_2\ga^ti_2i_1$.
Since $V(2)$ and $V(3)$ are  $V(0)$-module spectra, $\th(\be)=0$ and $\th(\ga)=0$ by \cite[Lemma 2.3]{t:zp}.
Similarly, $\th(i_k)=0$ and $\th(j_k)=0$ for $k=1,2$.
Therefore, \cite[Lemma 2.2]{t:zp} implies $\th(\be_{(s)})=0$ and $\th(\ga_{(t)})=0$. Therefore, $\be_{(s)}\ga_{(t)}=\ga_{(t)}\be_{(s)}$ by \cite[Cor. 2.7]{t:zp} as desired.

\subsection{On the action of $\ga_1$}

Note that $\ga_1=\al_1\be_{p-1}$.
Then, $\al_1\ga_1=\al_1^2\be_{p-1}=0$, $\lrk{\al_1,\al_1,\bes^p}\be_1\ga_1=-\al_1\lrk{\al_1,\al_1,\bes^p}\be_1\be_{p-1}=-\lrk{\al_1,\al_1,\al_1}\bes^p\be_1\be_{p-1}=0$ since $\lrk{\al_1,\al_1,\al_1}=0$, and
 $\lrk{\ga_1,p,\be_1}=\be_{p-1}\lrk{\al_1,p,\be_1}=\be_{p-1}j\underline{\al j_1}\be i_1i=0$.

For $t\ge 2$,
\sk{
	\be_{t}\en \de_{(1,1),1}\de_{(1,1),2}(v_2^{t})\es \de_{(1,1),1}([tv_2^{t-1}t_1^p+\C t2v_1v_2^{t-2}t_1^{2p}+v_1^2x])\\
	\cn [t(t-1)v_2^{t-2}t_2\ox t_1^p-tv_2^{t-1}b_0+\C t2v_2^{t-2}t_1\ox t_1^{2p}]\mod (p,v_1)\\
	\cn t(t-1)v_2^{t-2}k_0-tv_2^{t-1}b_0\mod (p,v_1)
}
and $\al_1\be_2\be_{p-1}\in E_2^5(S^0)$ is projected to $h_0(2k_0-2v_2b_0)(2v_2^{p-3}k_0+v_2^{p-2}b_0)=-2v_2^{p-2}h_0k_0b_0-2h_0v_2^{p-1}b_0^2$ in $E_2^5(V(2))$ under the induced map $i_*$ from the inclusion $i\cln S^0\to V(2)$ to the bottom cell.
Here, $k_0=[t_2\ox t_1^p+\frac12 t_1\ox t_1^{2p}]$.
Then, this element is detected by $-2v_2^{p-2}k_0\in E_1^3=E_2^{2,(p^2+p-1)q}(X\sm V(2))$ in the small descent \ss.
The killer of this element, if any, stays in the $E_1$-terms $E_1^2=E_2^{2,(p^2+p)q}(X\sm V(2))$, $E_1^1=E_2^{3,(p^2+2p-1)q}(X\sm V(2))$ and $E_1^0=E_2^{4,(p^2+2p)q}(X\sm V(2))$.
These are zero, and we see that the product is not zero.

%%%%%%%%%%%%%%%%%%%%   End of main body of article
%
%                             References
%
%   BiBTeX users uncomment the following line:
%
%\bibliographystyle{gtart}
%

\bibliographystyle{amsplain}

\end{document}